\documentclass[12pt,reqno]{amsart}

\usepackage{setspace, hyperref} \usepackage[dvips]{graphicx}

\newcommand\N {{\mathbb N}} 

\newcommand\R {\mathbb R}

\newcommand\astr{{}^\ast\hspace*{-0.7pt}\R}

\newcommand{\hr} {{}^{\mathfrak{h}}\hspace*{-0.6pt}\R}

\newcommand{\st}{\textbf{st}}

\DeclareMathOperator{\adequal}{\;\raisebox{-3pt}{$\ulcorner\!\urcorner$}\;}

\author[T.B.]{Tiziana Bascelli} \address{T. Bascelli, Lyceum Gymnasium
``F. Corradini'', Thiene, Italy}\email{tiziana.bascelli@virgilio.it}

\author[P.B.]{Piotr B\l{}aszczyk}\address{P. B\l{}aszczyk, Institute
of Mathematics, Pedagogical University of
Cracow,Poland}\email{pb@up.krakow.pl}

\author[V.K.]{Vladimir Kanovei} \address{V. Kanovei, IPPI, Moscow,
and MIIT, Moscow, Russia}\email{kanovei@googlemail.com}

\author[K.K.]{Karin U. Katz}\address{K. Katz, Department of
Mathematics, Bar Ilan University, Ramat Gan 52900
Israel}\email{katzmik@math.biu.ac.il}

\author[M.K.]{Mikhail G. Katz}\address{M. Katz, Department of
Mathematics, Bar Ilan University, Ramat Gan 52900
Israel}\email{katzmik@macs.biu.ac.il}

\author[S.K.]{Semen S. Kutateladze}\address{S. Kutateladze, Sobolev
Institute of Mathematics, Novosibirsk State University, Russia}
\email{sskut@math.nsc.ru}

\author[T.N.]{Tahl Nowik}\address{T. Nowik, Department of
Mathematics, Bar Ilan University, Ramat Gan 52900
Israel}\email{tahl@math.biu.ac.il}

\author[D.Sc.]{David M. Schaps}\address{D. Schaps, Department of
Classical Studies, Bar Ilan University, Ramat Gan 5290002
Israel}\email{dschaps@mail.biu.ac.il}

\author[D.Sh.]{David Sherry}\address{D. Sherry, Department of
Philosophy, Northern Arizona University, Flagstaff, AZ 86011,
US}\email{David.Sherry@nau.edu}

\numberwithin{equation}{section}

\begin{document}


\title{Gregory's sixth operation}

\begin{abstract}
In relation to a thesis put forward by Marx Wartofsky, we seek to show
that a historiography of mathematics requires an analysis of the
ontology of the part of mathematics under scrutiny.  Following Ian
Hacking, we point out that in the history of mathematics the amount of
contingency is larger than is usually thought.  As a case study, we
analyze the historians' approach to interpreting James Gregory's
expression \emph{ultimate terms} in his paper attempting to prove the
irrationality of~$\pi$.  Here Gregory referred to the \emph{last} or
\emph{ultimate terms} of a series.  More broadly, we analyze the
following questions: which modern framework is more appropriate for
interpreting the procedures at work in texts from the early history of
infinitesimal analysis?  as well as the related question: what is a
logical theory that is close to something early modern mathematicians
could have used when studying infinite series and quadrature problems?
We argue that what has been routinely viewed from the viewpoint of
classical analysis as an example of an ``unrigorous" practice, in fact
finds close procedural proxies in modern infinitesimal theories.  We
analyze a mix of social and religious reasons that had led to the
suppression of both the religious order of Gregory's teacher degli
Angeli, and Gregory's books at Venice, in the late 1660s.

\end{abstract}

\maketitle
\tableofcontents

\section{Introduction}
\label{s1}

Marx Wartofsky pointed out in his programmatic contribution \emph{The
Relation between Philosophy of Science and History of Science} that
there are many distinct possible relations between philosophy of
science and history of science, some ``more agreeable'' and fruitful
than others \cite[p.\;719ff]{Wa76}.  Accordingly, a fruitful relation
between history and philosophy of science requires a rich and complex
\emph{ontology} of that science.  In the case of mathematics, this
means that a fruitful relation between history and philosophy must go
beyond offering an ontology of the domain over which a certain piece
of mathematics ranges (say, numbers, functions, sets, infinitesimals,
structures, etc.).  Namely, it must develop the ontology of
mathematics \emph{as a scientific theory} itself (ibid., p.\;723).  A
crucial distinction here is that between the (historically relative)
\emph{ontology} of the mathematical objects in a certain historical
setting, and its \emph{procedures}, particularly emphasizing the
different roles these components play in the history of mathematics.
More precisely, \emph{procedures} serve as a representative of what
Wartofsky called the \emph{praxis} characteristic of the mathematics
of a certain time period, and \emph{ontology} in the narrow sense
takes care of the mathematical entities recognized at that time.  On
the procedure/entity distinction, A.\;Robinson had this to say:
\begin{quote}
\ldots from a formalist point of view we may look at our theory
syntactically and may consider that what we have done is to introduce
\emph{new deductive procedures} rather than new mathematical entities.
\cite[p.\;282]{Ro66} (emphasis in the original)
\end{quote}
As a case study, we analyze the text \emph{Vera Circuli} \cite{Gr67}
by James Gregory.

\section{Ultimate terms and termination of series}
\label{honor1}

Gregory studied under Italian indivisibilists%
\footnote{Today scholars distinguish carefully between indivisibles
(i.e., codimension one objects) and infinitesimals (i.e., of the same
dimension as the entity they make up); see e.g., \cite{Ko54}.
However, in the 17th century the situation was less clearcut.  The
term \emph{infinitesimal} itself was not coined until the 1670s; see
\cite{KS1}.}
and specifically Stefano degli Angeli during his years 1664--1668 in
Padua.  Some of Gregory's first books were published in Italy.  He
mathematical accomplishments include the series expansions not only
for the sine but also for the tangent and secant functions
\cite{Go11}.

The \emph{Vera Circuli} contains a characterisation of the
``termination'' of a convergent \emph{series} (i.e., \emph{sequence}
in modern terminology).  This was given by Gregory in the context of a
discussion of a double sequence (lower and upper bounds) of successive
polygonal approximations to the area of a circle:
\begin{quote}
\& igitur imaginando hanc seriem in infinitum continuari, possimus
imaginari vltimos terminos couergentes [sic] esse equales, quos
terminos equales appellamus seriei
terminationem. \cite[p.\;18-19]{Gr67}
\end{quote}
In the passage above, Gregory's \emph{seriem} refers to a
\emph{sequence}, and the expression \emph{terminus} has its usual
meaning of a \emph{term} of a sequence.  The passage can be rendered
in English as follows:
\begin{quote}
And so by imagining this series [i.e., sequence] to be continued to
infinity, we can imagine the ultimate convergent terms \emph{to be
equal}; and we call those equal ultimate terms the termination of the
series. [emphasis added]
\end{quote}
\cite[p.\;225]{Lu14} denotes the lower and upper bounds respectively
by~$I_n$ (for \emph{inscribed}) and~$C_n$ (for \emph{circumscribed}).
Gregory proves the recursive formulas~$I_{n+1}^2=C_n I_n$ and
$C_{n+1}=\frac{2C_n I_{n+1}}{C_n+I_{n+1}}$.  Gregory states that the
``ultimate convergent terms'' of the sequences~$I_n$ and~$C_n$ are
\emph{equal}.

After having defined the two series of inscribed and circumscribed
polygons, Gregory notes: 
\begin{quote}
atque in infinitum illam [=hanc polygonorum seriem] continuando,
manifestum est tandem exhiberi quantitatem sectori circulari,
elliptico vel hyperbolico ABEIOP
\ae quale[m]; differentia enim polygonorum complicatorum in seriei
continuatione semper diminuitur, it\`a vt omni exhibita quantitate
fieri possit minor, \& vt in sequenti theorematis Scholio
demonstrabimus: si igitur
pr\ae dicta polygonorum series terminari posset, hoc est, si
inueniretur vltimum illud polygonum inscriptum (si it\`a loqu\`i
liceat)
\ae quale vltimo illi polygono circumscripto, daretur infallibiliter
circuli \& hyperbol\ae{} quadra\-tura: sed quoniam difficile est, \&
in geometria omnin\`o fortasse inauditu[m] tales series terminare;
pr\ae mitte[n]d\ae{} sunt qu\ae{} dam propositiones \`e quibus
inueniri possit huiusmodi aliquot serierum terminationes, \& tandem
(si fieri possit) generalis methodus inueniendi omnium serierum
co[n]uergentium terminationes.
\end{quote}
This can be translated as follows: 
\begin{quote}
and that [series of polygons] being continued to infinity, it is clear
that a quantity equal to a circular, elliptic, or hyperbolic sector
ABEIOP will be produced.  The difference between [two $n$-th terms] in
the continuation of the series of complicated polygons always
diminishes so that it can become less than any given quantity indeed,
as we will prove in the Scholium to the theorem.  Thus, if the
abovementioned series of polygons can be terminated, that is, if that
ultimate inscribed polygon is found to be equal (so to speak) to that
ultimate circumscribed polygon, it would undoubtedly provide the
quadrature of a circle as well as a hyperbola.  But since it is
difficult, and in geometry perhaps unheard-of, for such a series to
come to an end [lit.: be terminated], we have to start by showing some
Propositions by means of which it is possible to find the terminations
of a certain number of series of this type, and finally (if it can be
done) a general method of finding terminations of all convergent
series.
\end{quote}
The passage clearly shows that Gregory is using the term ``ultimate
(or last) circumscribed polygon'' in a figurative sense, as indicated
by
\begin{itemize}
\item
his parenthetical `so to speak,' which indicates that he is not using
the term literally;
\item
his insistence that ``in geometry it is unheard-of'' for a sequence to
come to be terminated.
\end{itemize}
He makes it clear that he is using the word `termination' in a new
sense, which is precisely his sixth operation, as discussed below.

One possible interpretation of \emph{ultimate terms} would be the
following.  This could refer to those terms that are all closer than
epsilon to one another.  If ordinary terms are \emph{further} than
epsilon, that would make them different.  The difficulty for this
interpretation is that, even if ordinary terms are \emph{closer} than
epsilon, they are still \emph{different}, contrary to what Gregory
wrote about their being \emph{equal}.  M. Dehn and E.\;Hellinger
attribute to Gregory
\begin{quote}
a very general, new analytic process which he coordinates as the
``sixth'' operation along with the five traditional operations
(addition, subtraction, multiplication, division, and extraction of
roots).  In the introduction, he proudly states ``ut hae c nostra
inventio addat arithmeticae aliam operationem et geometriae aliam
rationis speciem, ante incognitam orbi geometrico.''  This operation
is, as a matter of fact, our modern limiting process.
\cite[p.\;157--158]{DH}
\end{quote}
We will have more to say about what this sixth operation could be
\emph{as a matter of fact} (see Section~\ref{three} on
\emph{shadow}-taking).  A. Malet expressed an appreciation of
Gregory's contribution to analysis in the following terms:
\begin{quote}
Studying Gregorie's work on ``Taylor'' expansions and his analytical
method of tangents, which has passed unnoticed so far, [we argue] that
Gregorie's work is a counter-example to the standard thesis that
geometry and algebra were opposed forces in 17th-century
mathematics. \cite[p.\;1]{Ma89}
\end{quote}
What is, then, Gregory's \emph{sixth operation} mentioned by Dehn and
Hellin\-ger, and how is it related to convergence?

\section{Law of continuity}
\label{two}

The use of infinity was not unusual for this period.  As we mentioned
in the introduction, Gregory fit naturally in the proud Italian
tradition of the method of indivisibles, and was a student of Stefano
degli Angeli at Padua between 1664 and 1668.  Degli Angeli published
sharp responses to critiques of indivisibles penned by jesuits Mario
Bettini and Andr\'e Tacquet.  Bettini's criticisms were extensions of
earlier criticisms by jesuit Paul Guldin.  Degli Angeli defended the
method of indivisibles against their criticisms.

Both indivisibles and degli Angeli himself appear to have been
controversial at the time in the eyes of the jesuit order, which
banned indivisibles from being taught in their colleges on several
occassions.  Thus, in 1632 (the year Galileo was summoned to stand
trial over heliocentrism) the Society's Revisors General led by Jacob
Bidermann banned teaching indivisibles in their colleges \cite{Fe90},
\cite[p.\;198]{Fe92}.  Indivisibles were placed on the Society's list
of \emph{permanently} banned doctrines in 1651 \cite{He96}.

It seems that Gregory's 1668 departure from Padua was well timed, for
his teacher degli Angeli's jesuat order%
\footnote{This was an older order than the jesuits.  Cavalieri had
also belonged to the jesuat order.}
was suppressed by papal brief in the same year, cutting short degli
Angeli's output on indivisibles.  Gregory's own books were suppressed
at Venice, according to a letter from John Collins to Gregory dated 25
november 1669, in which he writes:
\begin{quote}
One Mr.\;Norris a Master's Mate recently come from Venice, saith it
was there reported that your bookes were suppressed, not a booke of
them to be had anywhere, but from Dr.\;Caddenhead to whom application
being made for one of them, he presently sent him one (though a
stranger) refusing any thing for it. \cite[p.\;74]{Tu39}
\end{quote}
In a 1670 letter to Collins, Gregory writes:
\begin{quote}
I shall be very willing ye writ to Dr Caddenhead in Padua, for some of
my books.  In the mean time, I desire you to present my service to
him, and to inquire of him if my books be suppressed, and the reason
thereof.  (Gregory to Collins, St Andrews, March 7, 1670, in Turnbull
p.\;88)
\end{quote}
In a letter to Gregory, written in London, 29 september 1670, Collins
reported as follows: ``Father Bertet%
\footnote{\label{f13}Jean Bertet (1622-1692), jesuit, quit the Order
in\;1681.  In 1689 Bertet conspired with Leibniz and Antonio
Baldigiani in Rome to have the ban on Copernicanism
lifted. \cite{Wa12}}
sayth your Bookes are in great esteeme, but not to be procured in
Italy.'' (Turnbull p.\;107)

The publishers' apparent reluctance to get involved with Gregory's
books may also explain degli Angeli's silence on indivisibles
following the suppression of his order, but it is hard to say anything
definite in the matter until the archives at the Vatican dealing with
the suppression of the jesuats are opened to independent researchers.
Certainly one can understand Gregory's own caution in matters
infinitesimal (of course, the latter term wasn't coined until later).

John Wallis introduced the symbol~$\infty$ for an infinite number in
his \emph{Arithmetica Infinitorum} \cite{Wa56} and exploited an
infinitesimal number of the form~$\frac{1}{\infty}$ in area
calculations \cite[p.\;18]{Sc81}, over a decade before the publication
of Gregory's \emph{Vera Circuli}.  At about the same time, Isaac
Barrow ``dared to explore the logical underpinnings of
infinitesimals,'' as Malet put it:
\begin{quote}
Barrow, who dared to explore the logical underpinnings of
infinitesimals, was certainly modern and innovative when he publicly
defended the new mathematical methods against Tacquet and other
mathematical ``classicists'' reluctant to abandon the Aristotelian
continuum.  And after all, to use historical hindsight, it was the
non-Archimedean structure of the continuum linked to the notion of
infinitesimal and advocated by Barrow that was to prove immensely
fruitful as the basis for the Leibnizian differential
calculus. \cite[p.\;244]{Ma89}.
\end{quote}

We know that G. W. Leibniz was an avid reader of Gregory; see e.g.,
\cite{Le72}.  To elaborate on the link to Leibniz mentioned by Malet,
note that Leibniz might have interpreted Gregory's definition of
convergence as follows.  Leibniz's \emph{law of continuity}
\cite[p.\;93--94]{Le02a} asserts that whatever succeeds in the finite,
succeeds also in the infinite, and vice versa; see \cite{KS1} for
details.  Thus, if one can take terms of a sequence corresponding to a
finite value of the index~$n$, one should also be able to take terms
corresponding to infinite values of the index~$n$.  What Gregory
refers to as the ``ultimate'' terms would then be the terms~$I_n$
and~$C_n$ corresponding to an infinite index~$n$.

Leibniz interpreted equality as a relation in a larger sense of
equality \emph{up to} (negligible terms).  This was codified as his
\emph{transcendental law of homogeneity} \cite{Le10b}; see
\cite[p.\;33]{Bos} for a thorough discussion.  Thus, Leibniz wrote:
\begin{quote}
Caeterum aequalia esse puto, non tantum quorum differentia est omnino
nulla, sed et quorum differentia est incomparabiliter parva; et licet
ea Nihil omnino dici non debeat, non tamen est quantitas comparabilis
cum ipsis, quorum est differentia.  \cite[p. 322]{Le95b}
\end{quote}
This can be translated as follows: 
\begin{quote}
``Furthermore I think that not only those things are equal whose
difference is absolutely zero, but also whose difference is
incomparably small.  And although this [difference] need not
absolutely be called Nothing, neither is it a quantity comparable to
those whose difference it~is.''
\end{quote}
In the 17th century, such a generalized notion of equality was by no
means unique to Leibniz.  Indeed, Leibniz himself cites an antecedent
in Pierre de Fermat's technique (known as the method of
\emph{adequality}; see \cite{13e}), in the following terms:
\begin{quote}
Quod autem in aequationibus Fermatianis abjiciuntur termini, quos
ingrediuntur talia quadrata vel rectangula, non vero illi quos
ingrediuntur simplices lineae infinitesimae, ejus ratio non est quod
hae sint aliquid, illae vero sint nihil, sed quod termini ordinarii
per se destruuntur.%
\footnote{Translation: ``But the fact that in Fermat's equations those
terms into which such things enter as squares or rectangles [i.e.,
multiplied by themselves or by each other] are eliminated but not
those into which simple infinitesimal lines [i.e., segments]
enter--the reason for that is not because the latter are something
whereas the former are really nothing [as Nieuwentijt maintained], but
because ordinary terms cancel each other out.''}
\cite[p.\;323]{Le95b}
\end{quote}
On this page, Leibniz describes Fermat's method in a way similar to
Leibniz's own.  On occasion Leibniz used the notation ``$\adequal$''
for the relation of equality.  Note that Leibniz also used our symbol
``$=$'' and other signs for equality, and did not distinguish between
``$=$'' and ``$\adequal$'' in this regard.  To emphasize the special
meaning \emph{equality} had for Leibniz, it may be helpful to use the
symbol $\adequal$ so as to distinguish Leibniz's equality from the
modern notion of equality ``on the nose.''  Then Gregory's comment
about the equality of the ultimate terms translates into
\begin{equation}
\label{21a}
I_n\adequal{}C_n
\end{equation}
when~$n$ is infinite.

From the viewpoint of the modern Weierstrassian framework, it is
difficult to relate to Gregory's insight.  Thus, G. Ferraro translates
Gregory's ``vltimos terminos conuergentes'' as ``last convergent
terms'' \cite[p.\;21]{Fe08}, and goes on a few pages later to mention
\begin{quote}
Gregory's reference to \emph{the last term}, p.\;21. \ldots{} In
Leibniz they appear in a clearer way.  \cite[p.\;27, note~41]{Fe08}
(emphasis added)
\end{quote}
Ferraro may have provided an accurate translation of Gregory's
comment, but Ferraro's assumption that there is something unclear
about Gregory's comment because of an alleged ``last term'', is
unjustified.  Note that Ferraro's use of the singular ``last term''
(note~41) is not consistent with Gregory's use of the plural
\emph{terminos} (terms) in his book.  One may find it odd for a
mathematician of Gregory's caliber to hold that there is literally a
\emph{last} term in a sequence.  Dehn and Hellinger mention only the
plural ``last convergent terms'' \cite[p.\;158]{DH}.

\section{The Unguru controversy}
\label{three}

There is a debate in the community of historians whether it is
appropriate to use modern theories and/or modern notation in
interpreting mathematical texts of the past, with S. Unguru a staunch
opponent, whether with regard to interpreting Euclid, Apollonius, or
Fermat \cite{Un76}.  See \cite{Co13} for a summary of the debate.
Note that Ferraro does not follow Unguru in this respect.  Indeed,
Ferraro exploits the modern notation
\begin{equation}
\label{21}
\sum_{i=1}^\infty a_i
\end{equation}
for the sum of the series, already on page 5 of his book, while
discussing late 16th (!) century texts of Vi\`ete.  We note the
following two aspects of the notation~\eqref{21}:
\begin{enumerate}
\item
\label{i1}
It presupposes the modern epsilontic notion of limit, where
$S=\sum_{i=1}^\infty a_i$ means~$\forall\epsilon>0\;\exists
N\!\in\N\left(n>N\Longrightarrow\left|S-\sum_{i=1}^{n\phantom{I}}
a_i\right|<\epsilon\right)$, in the context of a Weierstrassian
framework involving a strictly Archimedean punctiform continuum;
\item
The symbol ``$\infty$'' occurring in Ferraro's formula has no meaning
other than a reminder that a limit was taken in the construction.  In
particular, this usage of the symbol~$\infty$ is distinct from its
original 17th century usage by Wallis, who used it to denote a
specific infinite number, and proceeded to work with infinitesimal
numbers like~$\frac{1}{\infty}$ (see Section~\ref{two}).
\end{enumerate}

We will avoid choosing sides in the debate over Unguru's proposal.%
\footnote{The sources of such a proposal go back (at least) to
A.\;Koyr\'e who wrote: ``Le probl\`eme du langage \`a adopter pour
l'exposition des oeuvres du pass\'e est extr\^emement grave et ne
comporte pas de solution parfaite.  En effet, si nous gardons la
langue (la terminologie) de l'auteur \'etudi\'e, nous risquons de le
laisser incompr\'ehensible, et si nous lui substituons la n\^otre, de
le trahir.''  \cite[p.\;335, note\;3]{Ko54}.}
However, once one resolves to exploit modern frameworks involving
punctiform continua/number systems, as Ferraro does, to interpret 17th
century texts, one still needs to address the following important
question:
\begin{quote}
\emph{Which modern framework is more appropriate for interpreting the
said historical texts?}
\end{quote}
Here appropriateness could be gauged in terms of providing the best
proxies for the \emph{procedural} moves found in the great 17th
century masters.  

\cite{Ha14} points out that there is a greater amount of contingency
in the historical evolution of mathematics than is generally thought.
Hacking proposes a \emph{Latin model} of development (of a natural
language like Latin, with the attendant contingencies of development
due to social factors) to the usual \emph{butterfly model} of
development (of a biological organism like a butterfly, which is
genetically predetermined inspite of apparently discontinuous changes
in its development).  This tends to undercut the apparent
inevitability of the Weierstrassian model.

We leave aside the \emph{ontological} or foundational questions of how
to justify the entities like points or numbers (in terms of modern
mathematical foundations), and focus instead of the \emph{procedures}
of the historical masters, as discussed in Section~\ref{s1}.

More specifically, is a modern Weierstrassian framework based on an
Archimedean continuum more appropriate for interpreting their
procedures, or is a modern infinitesimal system more appropriate for
this purpose?

Note that in a modern infinitesimal framework such as Robinson's,
sequences possess terms with infinite indices.  Gregory's relation can
be formalized in terms of the standard part principle in Robinson's
framework \cite{Ro66}.  This principle asserts that every finite
hyperreal number is infinitely close to a unique real number.

In more detail, in a hyperreal extension $\R\hookrightarrow\astr$ one
considers the set~$\hr\subseteq\astr$ of \emph{finite} hyperreals.
The set $\hr$ is the domain of the standard part function (also called
the \emph{shadow}) $\st\colon\hr\to\R$ rounding off each finite
hyperreal number to its nearest real number.

In the world of James Gregory, if each available term with an infinite
index~$n$ is indistinguishable (in the sense of being infinitely
close) from some standard number, then we ``terminate the series" (to
exploit Gregory's terminology) with this number, meaning that this
number is the limit of the sequence.  Gregory's definition corresponds
to a relation of infinite proximity in a hyperreal framework.  Namely
we have
\begin{equation}
\label{32}
I_n \approx C_n,
\end{equation}
where~$\approx$ is the relation of being infinitely close (i.e., the
difference is infinitesimal), and the common standard part of these
values is the limit of the sequence.
Equivalently,~$\st(I_n)=\st(C_n)$.  Mathematically speaking, this is
equivalent to a Weierstrassian epsilontic paraphrase along the lines
of item~\eqref{i1} above.

Recently Robinson's framework has become more visible thanks to
high-profile advocates like Terry Tao; see e.g., his work \cite{Ta14},
\cite{TV}.  The field has also had its share of high-profile
detractors like Errett Bishop and Alain Connes.  Their critiques were
critically analyzed in \cite{KK11d}, \cite{KL13}, and \cite{KKM}.
Further criticisms by J. Earman, K. Easwaran, H. M. Edwards, Ferraro,
J. Gray, P. Halmos, H. Ishiguro, G. Schubring, and Y. Sergeyev were
dealt with respectively in the following recent texts: 
\cite{KS1}, \;
\cite{Ba14}, \; 
\cite{KKKS}, \; 
\cite{Ba17}, \, 
\cite{Bl16c}, \,
\cite{Bl16a}, 
\cite{Ba16b},
\cite{Bl16d}, 
\cite{Gu16}.
In \cite{BK} we analyze the Cauchy scholarship of Judith Grabiner.
For a fresh look at Simon Stevin see \cite{KK12}.

\section{Conclusion}

We note a close fit between Gregory's procedure~\eqref{21a} and
procedure~\eqref{32} available in a modern infinitesimal framework.
The claim that ``[Gregory's] definition is rather different from the
modern one'' \cite[p.\;20]{Fe08} is only true with regard to a
\emph{Weierstrassian} modern definition.  Exploiting the richer syntax
available in a modern infinitesimal framework where Gregory's
procedure acquires a fitting proxy, it is possible to avoid the
pitfalls of attributing to a mathematician of Gregory's caliber odd
beliefs in an alleged ``last'' term in a sequence.

An infinitesimal framework also enables an interpretation of the
notion of ``ultimate terms'' as proxified by terms with infinite
index, and ``termination of the series'' as referring to the
assignable number infinitely close to a term with an infinite index,
by Leibniz's transcendental law of homogeneity (or the standard part
principle of Robinson's framework).

While some scholars seek to interpret Gregory's procedures in a
default modern post-Weierstrassian framework, arguably a modern
infinitesimal framework provides better proxies for Gregory's
procedural moves than a modern Weierstrassian one.

\section*{Acknowledgments}
M. Katz was partially supported by the Israel Science Foundation grant
no.\;1517/12.

\bigskip

\textbf{Tiziana Bascelli} graduated in mathematics (1993) and
philosophy (2006), and obtained a PhD in theoretical philosophy (2010)
from University of Padua, Italy.  She co-authored \emph{Campano da
Novara's Equatorium Planetarum.  Transcription, Italian translation
and commentary} (Padua, 2007) and \emph{Gali\-leo's `Sidereus Nuncius'
or `A Sidereal Message'}, translated from the Latin by William
R. Shea, introduction and notes by William R. Shea and Tiziana
Bascelli.  Sagamore Beach, MA (USA): Science History Publications,
2009.  She is an independent researcher in integrated history and
philosophy of science, in the field of early modern mechanics and
mathematics.  Her research interests are in the development of
infinitesimal objects and procedures in Seventeenth-century
mathematics.

\medskip

\textbf{Piotr B\l aszczyk} is Professor at the Institute of
Mathematics, Pedagogical University (Cracow, Poland).  He obtained
degrees in mathematics (1986) and philosophy (1994) from Jagiellonian
University (Cracow, Poland), and a PhD in ontology (2002) from
Jagiellonian University.  He authored \emph{Philosophical Analysis of
Richard Dedekind's memoir \emph{Stetigkeit und irrationale Zahlen}}
(2008, Habilitationsschrift).  He co-authored \emph{Euclid,
\emph{Elements, Books V--VI}.  Translation and commentary}, 2013; and
\emph{Descartes, Geometry.  Translation and commentary} (Cracow,
2015).  His research interest is in the idea of continuum and
continuity from Euclid to modern times.

\medskip

\textbf{Vladimir Kanovei} graduated in 1973 from Moscow State
University, and obtained a Ph.D. in physics and mathematics from
Moscow State University in 1976. In 1986, he became Doctor of Science
in physics and mathematics at Moscow Steklov Mathematical Institute
(MIAN).  He is currently Principal Researcher at the Institute for
Information Transmission Problems (IITP) and Professor at Moscow State
University of Railway Engineering (MIIT), Moscow, Russia.  Among his
publications is the book \emph{Borel equivalence relations. Structure
and classification}. University Lecture Series 44. American
Mathematical Society, Providence, RI, 2008.

\medskip

\textbf{Karin U. Katz} (B.A. Bryn Mawr College, '82); Ph.D. Indiana
University, '91) teaches mathematics at Bar Ilan University, Ramat
Gan, Israel.  Among her publications is the joint article ``Proofs and
retributions, or: why Sarah can't \emph{take} limits'' published in
\emph{Foundations of Science}.

\medskip

\textbf{Mikhail G. Katz} (BA Harvard '80; PhD Columbia '84) is
Professor of Mathematics at Bar Ilan University, Ramat Gan, Israel.
He is interested in Riemannian geometry, infinitesimals, debunking
mathematical history written by the victors, as well as in true
infinitesimal differential geometry; see \emph{Journal of Logic and
Analysis} \textbf{7}:5 (2015), 1--44 at
\url{http://dx.doi.org/10.4115/jla.2015.7.5}

\medskip

\textbf{Semen S. Kutateladze} was born in 1945 in Leningrad (now
St.~Petersburg).  He is a senior principal officer of the Sobolev
Institute of Mathematics in Novosibirsk and professor at Novosibirsk
State University.  He authored more than 20 books and 200 papers in
functional analysis, convex geometry, optimization, and nonstandard
and Boolean valued analysis.  He is a member of the editorial boards
of \emph{Siberian Mathematical Journal}, \emph{Journal of Applied and
Industrial Mathematics}, \emph{Positivity}, \emph{Mathematical Notes},
etc.

\medskip

\textbf{Tahl Nowik} Ph.D. in Mathematics from the Hebrew University
1996.  Postdoctorate at Columbia University 1996--2000.  Member of the
Mathematics department at Bar Ilan University since 2000.  Research
interests are in low dimensional topology, finite type invariants,
stochastic topology, and nonstandard analysis.

\medskip

\textbf{David Schaps} is Professor Emeritus of Classical Studies at
Bar Ilan University, Israel.  His books include \emph{Economic Rights
of Women in Ancient Greece}, \emph{The Invention of Coinage and the
Monetization of Ancient Greece}, and \emph{Handbook for Classical
Research}; among his articles are ``The Woman Least Mentioned:
Etiquette and Women's Names", ``What was Free about a Free Athenian
Woman?'', and ``Zeus the Wife-Beater''.

\medskip

\textbf{David Sherry} is Professor of Philosophy at Northern Arizona
University, in the tall, cool pines of the Colorado Plateau.  He has
research interests in philosophy of mathematics, especially applied
mathematics and non-standard analysis.  Recent publications include
``Fields and the Intelligibility of Contact Action,'' \emph{Philosophy
90} (2015), 457--478.  ``Leibniz's Infinitesimals: Their Fictionality,
their Modern Implementations, and their Foes from Berkeley to Russell
and Beyond,'' with Mikhail Katz, \emph{Erkenntnis 78} (2013), 571-625.
``Infinitesimals, Imaginaries, Ideals, and Fictions,'' with Mikhail
Katz, \emph{Studia Leibnitiana 44} (2012), 166--192.  ``Thermoscopes,
Thermometers, and the Foundations of Measurement,'' \emph{Studies in
History and Philosophy of Science 24} (2011), 509--524.  ``Reason,
Habit, and Applied Mathematics,'' \emph{Hume Studies 35} (2009),
57-85.

\vfill\eject

\end{document}